\documentclass[12pt]{article}

\usepackage{epsfig}

\makeatletter
\def\eqalign#1{\null\vcenter{\def\\{\cr}\openup\jot\m@th
  \ialign{\strut$\displaystyle{##}$\hfil&$\displaystyle{{}##}$\hfil
      \crcr#1\crcr}}\,}
\makeatother

\newcommand{\be}{\begin{equation}} 
\newcommand{\ee}{\end{equation}}
\newcommand{\beq}{\begin{eqnarray}}
\newcommand{\eeq}{\end{eqnarray}}
\newcommand{\bt}{\begin{theorem}}
\newcommand{\et}{\end{theorem}}
\newcommand{\bl}{\begin{lemma}}
\newcommand{\el}{\end{lemma}}
\newcommand{\bc}{\begin{corollary}}
\newcommand{\ec}{\end{corollary}}
\newcommand{\ba}{\begin{array}}
\newcommand{\ea}{\end{array}}

\newcommand{\la}{\label}
\newcommand{\ci}{\cite}

\newtheorem{theorem}{Theorem}
\newtheorem{lemma}[theorem]{LEMMA}
\newtheorem{corollary}[theorem]{COROLLARY}

\newcommand{\De}{\Delta}
\newcommand{\al}{\alpha}
\newcommand{\ga}{\gamma}

\renewcommand{\th}{\theta}

\newcommand{\vp}{\varphi}

\newcommand{\bi}{\bibitem}

\newfont{\msbm}{msbm10 scaled\magstep1}
\newfont{\msbms}{msbm7 scaled\magstep1} 

\newcommand{\bbc}{\mbox{$\mbox{\msbm C}$}}

\begin{document}
\begin{center}
{\Large\bf Asymptotics for Toeplitz determinants on a circular arc}\\
\bigskip\bigskip
I. V. Krasovsky\\
\bigskip \bigskip
Technische Universit\"at Berlin, Germany\\
and\\
Brunel University West London, United Kingdom\\
\bigskip\bigskip
 
\end{center}
%

\section{Introduction}
In this work we find an asymptotic formula for a Toeplitz
determinant with the symbol supported on an arc of the unit circle in
the case when the symbol has Fisher-Hartwig singularities.  

Let $f(\theta)$ be an integrable function on the unit circle. It
was proved by Szeg\H o \ci{gs} that if $f(\th)$ is positive and
sufficiently smooth, namely, its derivative satisfies a Lipshitz
condition, (these requirements have been later relaxed)
then the Toeplitz determinant
\[
D_n(f)=\det\left({1\over 2\pi}
\int_0^{2\pi}e^{-i(j-k)\th}f(\th)d\th\right),
\qquad j,k=0,1\dots,n, 
\]
has the asymptotic expansion
\be
D_{n-1}(f)\sim H(f)^n E(f),\qquad n\to\infty,\la{S}
\ee
where
\be
\eqalign{
H(f)=\exp{(\ln f)_0},\qquad 
E(f)=\exp{\sum_{k=1}^\infty k(\ln f)_k (\ln f)_{-k}},\\
(\ln f)_k={1\over 2\pi}\int_0^{2\pi}e^{-ik\th}\ln f(\th)d\th,\qquad
k=0,1,\dots}
\ee
This formula is not valid for $f(\th)$ with zeros or singularities. 
The asymptotics in the case when  $f(\th)$ has zeros or root-type 
singularities were conjectured by Lenard \ci{l}, 
Fisher and Hartwig \ci{fh} and proved
by Widom \ci{Wsing}. Namely, let
\be
f(\th)=\psi(\th)\prod_{r=1}^R(2-2\cos(\th-\th_r))^{\al_r}=
\psi(\th)\prod_{r=1}^R\left|2\sin{\th-\th_r \over 2}\right|^{2\al_r}, 
\ee
where $\psi(\th)$ satisfies the conditions of Szeg\H o's theorem and all
$\al_r>-1/2$ (this is the most general case of real $\al_r$'s since 
it is a condition for existence of Fourier coefficients).
Then
\be
D_{n-1}(f)\sim H(\psi)^n n^{\sum_1^R\al_r^2}
E(\psi)\prod_1^R \left({H(\psi)\over
    \psi(\th_r)}\right)^{\al_r} {G(\al_r+1)^2\over G(2\al_r+1)}
 \prod_{r\neq s}\left|2\sin{\th_r-\th_s \over 2}
\right|^{-\al_r\al_s},\la{sing}
\ee
where $G(z)$ is Barnes' G-function.
This was further generalized for the case of $f(\th)$ with jumps by
Basor \ci{b}, B\"ottcher and Silbermann \ci{bs1}, and other authors.
An incomplete list includes the works \ci{b}--\ci{e}.
Note that the positivity condition for $f(\th)$ can be removed:
the asymptotics (\ref{sing}) are actually proved in \cite{Wsing}
for a complex-valued $f(\th)$ with $\psi(\th)\neq 0$ and 
the change of the argument over the closed circle
$\Delta_{-\pi\le\th\le\pi}\arg\psi(\th)=0$.  
 
It is interesting how the asymptotics look like in case when $f(\th)$
is zero on an arc. This question was addressed by Widom in \ci{Warc}.
Suppose that $f(\vp)$ is supported on an arc $\al\le\vp\le2\pi-\al$,
where it is positive, smooth  enough,
and symmetric: $f(\vp)=f(2\pi-\vp)$.
Then it is proved in \ci{Warc} that 
\be
D_{n-1}(f)\sim \ga^{n^2}H(F)^n\left(n\sin{\al\over
  2}\right)^{-1/4}{G(3/2)^2\over\sqrt{\pi}}E(F),\la{arc}
\ee
with $\ga=\cos(\al/2)$ and
\[
F(\th)=f(2\arccos(\ga\cos{\th\over 2})), 
\]
the function $f(\vp)$ ``extended'' to the whole circle. Note that as
is known
\[
{G(3/2)^2\over\sqrt{\pi}}= 2^{1/12}e^{3\zeta'(-1)},
\]
where $\zeta'(x)$ is the derivative of Riemann's zeta function.
Widom's proof was based on a theorem of Hirschman for a Toeplitz-like
determinant, where in the matrix elements 
the Fourier coefficients of $f(\th)$ on the circle, 
i.e. the coefficients w.r.t. the system $z^k$, are replaced by the 
coefficients of a function on the interval $[-1,1]$ w.r.t. the Legendre 
polynomials.

Just like (\ref{S}) this formula breaks down if $f(\vp)$ has zeros
or singularities on the arc. Resolution of this question is
presented here. Let
\be
f(\vp)=\psi(\vp)\prod_{r=1}^R(2-2\cos(\vp-\vp_r))^{\al_r}, \qquad
\al_r>-1/2,\quad\vp_1=\al,\quad \vp_R=2\pi-\al,\la{f} 
\ee
be supported on an arc $\al\le\vp\le2\pi-\al$, $f(\vp)=f(2\pi-\vp)$
(in particular, $\vp_r=2\pi-\vp_{R+1-r}$), 
and let $\psi(\vp)$ satisfy the conditions of Szeg\H o's theorem.
Then 
\be
\eqalign{
D_{n-1}(f)\sim \ga^{(n+\sum_1^R \al_r)^2}
H(\Psi)^n \left({n\over\ga}\right)^{2\al_1^2+\sum_1^R \al_r^2}
\left(n\sin{\al\over 2}\right)^{-1/4}
2^{4\al_1+2\al_1^2}\Gamma(1+2\al_1)
{G(3/2+2\al_1)^2\over\sqrt{\pi}G(2+4\al_1)}\times\\
E(\Psi)(\sin\al)^{2\al_1^2}
\prod_1^R\left({H(\Psi)\over
    \psi(\vp_r)}\right)^{\al_r}
\prod_{r=2\atop r\neq R}^{R-1}
 {G(\al_r+1)^2\over G(2\al_r+1)}
\left({\sin^2 (\vp_r/2) \over 1-\ga^{-2}\cos^2(\vp_r/2)}
\right)^{\al_r^2/2}
 \prod_{r\neq s}\left|2\sin{\vp_r-\vp_s \over 2}\right|^{-\al_r\al_s},}
\la{arcsing}
\ee
where
\[
\Psi(\th)=\psi(2\arccos(\ga\cos{\th\over 2})). 
\]

In particular, when zeros and singularities are absent, we reproduce
the asymptotics (\ref{arc}).

For technical reasons, we need to replace $\psi(\vp)$ by an analytic function.
Smoothness of $\psi(\vp)$ implies existence of a
trigonometric polynomial $\psi_0(\vp)$ such that 
$\mu(\vp)=\psi(\vp)-\psi_0(\vp)=O(1/C)$ uniformly in $\vp$ on the arc 
for some constant $C(n)$ which can be chosen arbitrary large.
Let, e.g., $C(n)=n^n$. Then a simple analysis of the r.h.s. of (\ref{arcsing})
(using smoothness and positivity of $\psi(\vp)$) and the l.h.s. 
of (\ref{arcsing}) (using Hadamard's inequality) shows that it is 
sufficient to prove the theorem when $\psi$ 
is replaced by $\psi_0$. To simplify notation, we just assume in what follows
that $\psi(\vp)$ is analytic. (Note that we have to bear in mind that 
$\psi(\vp)$ now depends on $n$ in a particular way, as it now stands for 
an analytic approximation to the original $\psi$.)

To obtain (\ref{arcsing}), we use (\ref{sing}) to evaluate 
$D_{n-1}(\hat F)$, where
\be 
\hat F(\th)=\frac{F(\th)\sin(\th/2)}{\sqrt{1-\ga^2\cos^2(\th/2)}},
\qquad F(\th)=f(2\arccos(\ga\cos{\th\over 2})),
\la{hF}
\ee
is a symbol on the
whole circle, and the following relation
\be
{D_{n-1}(f)\over D_{n-1}(\hat F)}=
\ga^{n^2+n}{P_n(\gamma^{-1})\over P_n(1)}
\sim
\ga^{n^2}\frac{\sqrt{\tilde f(\al)}(1+\sin(\al/2))^n}{\De(F)\sqrt{\pi
    n\sin(\al/2)}}
  \frac{2^{2\al_1}\Gamma(1+2\al_1)}{n^{2\al_1}}
\la{ratio},
\ee
where
\begin{eqnarray}
\tilde f(\vp)={(2\ga)^{4\al_1}f(\vp)\over 
\prod_{r=1,R}(2-2\cos(\vp-\vp_r))^{\al_r}},\la{tf}
\\
\ln\De(F)={\sin(\al/2) \over 4\pi}\int_0^{2\pi}\frac{\ln F(\th)}
{1-\ga^2\cos^2(\th/2)}d\th,\la{De1}
\end{eqnarray}
and $P_k(x)=x^k+\cdots$, 
$k=0,1,\dots$ are the monic polynomials orthogonal on $[-1,1]$
w.r.t. the weight function
\begin{eqnarray}
w(x)=\frac{f(2\arccos\ga x)}{\sqrt{1-\ga^2 x^2}}:\label{w}\\
\int_{-1}^1 P_n(x)P_m(x)w(x)dx=h_n\delta_{mn},\quad h_n>0,\quad
m,n=0,1,\dots
\end{eqnarray}

In our case of a symmetric weight ($w(x)=w(-x)$), these polynomials
are simply related (see \ci{szego,f,zhedanov,fredh}) to the polynomials
$\Phi_k(z,\beta)=z^k+\cdots$, $k=0,1,\dots$ orthogonal w.r.t. 
$g(\vp)=f(2\arccos(\ga\cos{\vp\over 2}/\cos{\beta\over 2}))$
on an arc:
\begin{equation}
{1\over 2\pi}\int_\beta^{2\pi-\beta} \Phi_n(z)\overline{\Phi_m(z)}
g(\vp)d\vp=
{\tilde h}_n\delta_{mn}; \quad z=e^{i\vp},\quad \tilde h_n>0,
\quad m,n=0,1,\dots\label{phi}
\end{equation}

The first equation in (\ref{ratio}) follows from the 
connection between 1) Toeplitz matrices, 2) orthogonal polynomials
on an arc of the unit circle, and 3) those on a segment of the real
axis.
To obtain the second equation, we need to analyze polynomials $P_k(z)$
asymptotically. This is partly achieved with the help of an approach
based on a matrix Riemann-Hilbert problem, a new method which turned
out to be
very effective for analysis of asymptotic problems \ci{fik,dz,d}.  
Most Riemann-Hilbert analysis we need is contained in
\ci{kuimcl,kuivan}
where the asymptotics were found for polynomials orthogonal on
$[-1,1]$ with a weight $(1-x)^a(1+x)^b h(x)$, $h(x)$ is positive and
real analytic. We get the main term in the asymptotics of $P_n(1)$ by
generalizing a small part of \ci{kuimcl,kuivan} to the weight (\ref{w}).
The main term for $P_n(\ga^{-1})$ can be obtained much easier, from the
well-known Szeg\H o asymptotics (\ci{szego}, Thm.12.1.2) which apply
here because $\ga^{-1}$ lies in the exterior of the orthogonality 
interval $[-1,1]$.
As will be clear from the proof, equations (\ref{ratio}) also hold
when $f(\th)$ has jumps in addition to zeros and root singularities
(see \ci{fh}). Therefore the result of Basor \ci{b} can also be
extended to symmetric symbols on an arc in the same way as (\ref{sing}).

\section{The ratio $D_{n-1}(f)/D_{n-1}(\hat F)$.}
In this section we obtain equations (\ref{ratio}).
Due to the above mentioned connection between polynomials 
orthogonal on an interval of
the real axis and on an arc of the unit circle 
one can write the following relations (Lemma 3.1. of \ci{fredh}):
\begin{eqnarray}
D_{n-1}(f)=2^{n(n-1)}\frac{\gamma^{n^2+n}}
{\pi^n}P_n(1/\gamma)
\prod_{j=0}^{n-1} h_j,\qquad f(\vp)=w(\ga^{-1}\cos{\vp\over2})
\sin{\vp\over2};
\label{D}\\
D_{n-1}(\hat F)=2^{n(n-1)}\frac{1}
{\pi^n}P_n(1)
\prod_{j=0}^{n-1} h_j,\qquad 
\hat F(\th)=w(\cos{\th\over2})\sin{\th\over2}.
\end{eqnarray}
In the first of them the corresponding arc is $\al\le\vp\le 2\pi-\al$,
and in the second one, the arc coincides with the whole circle
$0\le\th<2\pi$. Setting $\ga\cos(\th/2)=\cos(\vp/2)$ and comparing
the expressions for $f$ and $\hat F$ we get (\ref{hF}). 
The ratio of $D_{n-1}$'s yields the first equation in
(\ref{ratio}).  
Note that we can also rewrite it as follows (using, e.g., formula (14)
of \ci{fredh}):
\be
{D_{n-1}(f)\over D_{n-1}(\hat F)}=\frac{2 \Phi_n(1,\al)\ga^{n^2+n}}
{e^{-in\al/2}\Phi_n(e^{i\al},\al)+\mathrm{c.c.}}.
\ee
Thus, to obtain the second (asymptotic) equation in (\ref{ratio}),
we can either consider the asymptotics of polynomials on an arc or on
an interval. We choose the second option since in that case we can 
use the results of Kuijlaars, McLauphlin, Van Assche, and Vanlessen
\ci{kuimcl} instead of reformulating the 
Riemann-Hilbert solution for polynomials on an arc. 
As indicated in the introduction, we need these methods to obtain the
main term in the asymptotics of $P_n(1)$. 
In what follows, we mostly dwell on the differences of our case from 
\ci{kuimcl} and refer the reader to that work for other details.

Consider the $2\times 2$ matrix
\begin{equation} \label{RHPYsolution}
    Y(z) =
    \pmatrix{
   P_n(z) & \frac{1}{2\pi i} \int_{-1}^1 P_n(x) w(x){dx\over x-z} \cr
-2\pi i \kappa_{n-1}^2 P_{n-1}(z) & 
-\kappa_{n-1}^2 \int_{-1}^1 P_{n-1}(x)w(x){dx\over x-z}
    },
\end{equation}
where $\kappa_n=h_n^{-1/2}$ is the leading coefficient of the
orthonormal polynomials.
$Y(z)$ is the unique solution of the following Riemann-Hilbert problem
(the proof is almost the same as in \ci{kuimcl}, the difference
being additional singularities inside $(-1,1)$):
\begin{enumerate}
    \item[(a)]
        $Y(z)$ is  analytic for $z\in\bbc \setminus [-1,1]$.
    \item[(b)]
For $x \in (-1,1)\setminus\cup_r x_r$,  $x_r=\ga^{-1}\cos(\vp_r/2)$,
$Y$ has continuous boundary values
$Y_{+}(x)$ as $z$ approaches $x$ from
above, and $Y_{-}(x)$, from below. They are related by 
the jump condition
\begin{equation}\label{RHPYb}
            Y_+(x) = Y_-(x)
            \pmatrix{
                1 & w(x) \cr
                0 & 1},
            \qquad\mbox{$x \in (-1,1)\setminus\cup_r x_r$.}
        \end{equation}
    \item[(c)]
        $Y(z)$ has the following asymptotic behavior at infinity:
        \begin{equation} \label{RHPYc}
            Y(z) = \left(I+ O \left( \frac{1}{z} \right)\right)
            \pmatrix{
                z^{n} & 0 \cr
                0 & z^{-n}}, \qquad \mbox{as $z\to\infty$.}
        \end{equation}
    \item[(d)]
       Near the points $1=x_1$, $-1=x_{R}$, $x_r$ 
        \begin{equation}\label{RHPYd}
            Y(z) = \left\{
            \begin{array}{cl}
                O\pmatrix{
                    1 & |z-x_r|^{2\al_r} \cr 
                    1 & |z-x_r|^{2\al_r}},
                &\mbox{if $\al_r<0$,} \\[2ex]
                O\pmatrix{ 
                    1 & \log|z-x_r| \cr
                    1 & \log|z-x_r|},
                &\mbox{if $\al_r=0$,} \\[2ex]
                O\pmatrix{
                    1 & 1 \cr
                    1 & 1},
                &\mbox{if $\al_r>0$,}
            \end{array}\right.
        \end{equation}
as $z \to x_r$, $z \in \bbc \setminus [-1,1]$.
\end{enumerate}

Here the points $x_r$ correspond to $\vp_r$ under the mapping
$x=\ga^{-1}\cos(\vp/2)$. 
Moreover, recalling (\ref{f}), (\ref{w}), we have
\be
w(x)=\frac{\psi(2\arccos\ga x)}{\sqrt{1-\ga^2 x^2}}
\prod_1^R (2\ga)^{2\al_r}|x_r-x|^{2\al_r}.
\label{w2}
\ee

The proof of (\ref{RHPYd}) 
is similar to that in \ci{kuimcl}.

Now we find the asymptotics of the matrix element $Y_{11}(z)$ 
in the neighbourhood of
$z=1$. To do this, we apply a series of transformations to the 
Riemann-Hilbert problem for $Y(z)$, as is usual in this method.
These transformations $Y\mapsto T\mapsto S$ are the same as in
\ci{kuimcl} but the contour where the jump for $S(z)$ occurs is
different. It consists of $R-1$ lenses which touch at zeros and
singularities $x_r$ of the weight $w_0(x)$ (see the figure).
\begin{figure}
\centerline{\psfig{file=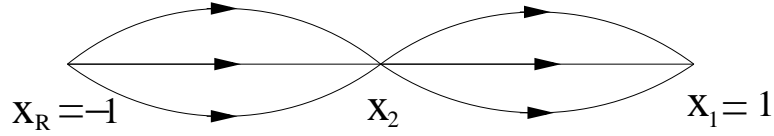,width=4.0in,angle=0}}
\vspace{0.1cm}
\caption{
The contour for the Riemann-Hilbert problem
is sketched in the case of  $R=3$.}
\label{fig1}
\end{figure}
In \ci{kuimcl} there was only one lens with the endpoints $-1$ and $1$.
We now consider the neighbourhoods of the points $x_r$ and the region
outside of them separately. In the outside region, a parametrix
$N(z)$ is given by the same expression as in \ci{kuimcl} where we
substitute our weight $w_0(x)$. It coincides with the main asymptotic 
term for $Y(z)$ obtained from Szeg\H o's Theorem 12.1.2 in \ci{szego}.
In the neighbourhood of
each point $x_r$, one could solve the local Riemann-Hilbert problem
such that it fits $N(z)$ on the neighbourhood's boundary. Of those
solutions, we need only the one in the neighbourhood of $x_1=1$. 
Indeed, it is easy to show again like in Theorems 7.8, 7.9 of 
\ci{deiftkriech} that corrections from neighbourhoods of the other
points $x_r$ contribute at most to the 
next-leading term in the asymptotics near $x=1$. Therefore
the leading term near $x=1$ coincides exactly with that for the
problem considered in \ci{kuimcl, kuivan}. It is given by the formula
(2.23) and the one after (3.7) of \ci{kuivan}. Taking there the limit
$x\to 1$, we obtain the asymptotics of $Y_{11}(1)$ and finally, recalling 
(\ref{RHPYsolution}), get the relation:
\be
P_n(1)\sim {{\cal D}_\infty \over 2^n}\sqrt{2\over 1+\sin(\al/2)}
\sqrt{\pi n\over v(1)}
{n^{2\al_1}\over 2^{2\al_1}\Gamma(1+2\al_1)},\la{P1}
\ee
where   
\be
\ln {\cal D}_\infty= {1\over 2\pi}\int_{-1}^1\frac{\ln f(2\arccos\ga x)}
{\sqrt{1-x^2}}dx=
{1\over 4\pi}\int_0^{2\pi}\ln F(\th)d\th\la{Dinf}
\ee
and $v(x)=w(x)/(1-x^2)^{2\al_1}$.

One can verify that there is no problem in the fact that
$\psi(\vp)$ (and hence $w(x)$) depends on $n$ (in a special way as 
it approximates some smooth function). 

We now turn to calculation of the asymptotics for $P_n(\ga^{-1})$.
Here it suffices to use an old result of Szeg\H o.
The leading term for asymptotics of {\it orthonormal} polynomials
$p_n(x)$ for $x$ outside $[-1,1]$ and a wide class of weights on
$[-1,1]$ (which includes our weight (\ref{w2})) is given by Theorem
12.1.2 in \ci{szego}. Namely,
\begin{eqnarray}
p_n(x)\sim {z^n \over \sqrt{2\pi}{\cal D}(z^{-1})},\qquad 
x={1\over 2}(z+z^{-1}),\la{p}\\ 
\ln {\cal D}(z)={1\over 4\pi}\int_{-\pi}^\pi \ln (w(\cos t)|\sin t|)
\frac{1+z e^{-it}}{1-z e^{-it}}dt.\la{Dz}
\end{eqnarray}  
For $x=1/\ga=1/\cos(\al/2)$ we have $z=(1+\sin(\al/2))/\ga\equiv z_0$.
Substituting (\ref{w}) into (\ref{Dz}), separating the integral into 
2 parts, one over $[0,\pi]$, the other over $[-\pi,0]$, and
changing the variables $\th\mapsto -\th$ in the second one, we get
\[
\ln {\cal D}(z_0^{-1})=
{1\over 4\pi}\int_0^\pi\ln\left(F(2\th)
{\sin\th \over\sqrt{1-\ga^2\cos^2\th}}\right) 
{2\sin(\al/2)\over {1-\ga\cos\th}}d\th=
\ln\De(F)+I,
\]
where
\[
\ln\De(F)={\sin(\al/2)\over 2\pi}\int_0^\pi
{\ln F(2\th)\over 1-\ga\cos\th}d\th=
{\sin(\al/2)\over 4\pi}\int_0^{2\pi}
{\ln F(\th)\over 1-\ga^2\cos^2(\th/2)}d\th
\]
(we used the symmetry $F(\th)=F(2\pi-\th)$
to get the last equation)
and
\be
I={\sin(\al/2)\over 4\pi}\int_0^{2\pi}
\ln{|\sin\th| \over\sqrt{1-\ga^2\cos^2\th}}
{d\th\over 1-\ga\cos\th}.\la{I}
\ee

The last integral can be evaluated explicitly.
We have for $|x|>1$, $0<a<1$
\be
-{\sqrt{x^2-1}\over 4\pi}\int_0^{2\pi}
{\ln(1-a\cos^2\th)\over x-\cos\th}d\th=\ln\frac{x+\sqrt{x^2-1}}
{\sqrt{x^2-1}+x\sqrt{1-a}}.\la{int0}
\ee
This integral is calculated by differentiating w.r.t. $a$ and then
computing residues in the complex plane of $e^{i\th}$.
Taking here $a\to 1$ and setting $x=1/\ga$ gives part of the integral
(\ref{I}) with $\ln|\sin\th|$. The part with $\ln(1-\ga^2\cos^2\th)$
is obtained by setting $a=\ga^2$, $x=1/\ga$. Thus we get 
$I=-\ln\sqrt{2}$, and therefore
\be
\ln {\cal D}(z_0^{-1})=\ln{\De(F)\over\sqrt{2}}.\la{Dres}
\ee
The main term in the asymptotics of the leading coefficient
$\kappa_n$ of
$p_n(x)=\kappa_n x^n+\cdots$ is given by Theorem 12.7.1 of \ci{szego}.
A calculation like the one just carried out yields
\be
\kappa_n\sim{2^n\over {\cal D}_\infty}\sqrt{1+\sin(\al/2)\over 2\pi},\la{kappa}
\ee
where ${\cal D}_\infty$ is given by (\ref{Dinf}).
Putting (\ref{p}), (\ref{Dres}), and (\ref{kappa}) together,
we finally obtain
\be
P_n(\ga^{-1})={p_n(\ga^{-1})\over\kappa_n}\sim
\left({1+\sin(\al/2) \over 2\ga}\right)^n 
\sqrt{2\over 1+\sin(\al/2)}{{\cal D}_\infty \over \De(F)}.\la{Pg}
\ee
Equations (\ref{P1}) and (\ref{Pg}) give 
the asymptotics (\ref{ratio}).

\section{Asymptotics for Toeplitz determinants}
We are now ready to obtain the asymptotics (\ref{arcsing}).
An important role will be played by the map between the arc and the
circle given by $\cos(\vp/2)=\ga\cos(\th/2)$.
Using the symmetry of the symbol $f(\vp)=f(2\pi-\vp)$, and hence
$\vp_r=2\pi-\vp_{R+1-r}$, we can write (\ref{f}) as ($\th_{r_0}=\pi$)
\be\eqalign{
F(\th)=f(\vp)=
\psi(\vp)\prod_{r=1}^{[R/2]}\left|4\sin{\vp-\vp_r\over2}
\sin{\vp+\vp_r\over2}\right|^{2\al_r}\left|2\cos{\vp\over2}
\right|^{2\al_{r_0}}=\\
\psi(\vp)\prod_{r=1}^{[R/2]}
\left(4\left|\cos^2{\vp\over2}-\cos^2{\vp_r\over2}\right|
\right)^{2\al_r}\left|2\cos{\vp\over2}
\right|^{2\al_{r_0}}=\\
\ga^{2\sum_r\al_r}\Psi(\th)\prod_{r=1}^R
\left|2\sin{\th-\th_r \over 2}\right|^{2\al_r}.}\la{fF}
\ee
Consider the function $\hat F$ (\ref{hF}).
Vanishing of $f(\vp)$ on the arc $-\al\le\vp\le\al$ is 
``replaced'' in $\th$ variable by a singularity of $\hat F$
at $\th=0$ given by
the factor $2\sin(\th/2)$. If $f(\vp)$ also has singularities at the
endpoints of the arc, the corresponding factor becomes  
$(2\sin(\th/2))^{1+4\al_1}$. 
We see that Widom's formula (\ref{sing}) holds for 
$\hat F$. The role of $\psi$ is now played by
\[
{\ga^{2\sum_r\al_r} \Psi(\th)
\over 2\sqrt{1-\ga^2\cos^2(\th/2)}},
\]
$\al_1$ in the formula is replaced by $1/2+2\al_1$, $R$, by $R-1$ (as
the ends of the arc merge in $\th$ variable).
Let us now simplify this expression for $D_{n-1}(\hat F)$.
We have
\be
\left(\ln2\sqrt{1-\ga^2\cos^2{\th\over2}}\right)_k=
\cases{-{1\over 2|k|}\left({1-\sin(\al/2)\over\ga}\right)^{2|k|}, 
        & $k\neq 0$\cr
        \ln(1+\sin(\al/2)),& $k=0$.}\la{int}
\ee
These integrals are calculated as (\ref{int0}).
In particular, if we set $\al=0$ for $k=0$, we get
\[
\left(\ln 2\left|\sin{\th-\th_r\over 2}\right|\right)_0=0.
\]

The first obvious application of these relations is the expression
\be
H(\hat F)=
H\left({\ga^{2\sum_r\al_r} \Psi(\th)
\over 2\sqrt{1-\ga^2\cos^2(\th/2)}}\right)=
{\ga^{2\sum_r\al_r}H(\Psi)\over 1+\sin(\al/2)}.
\ee
According to (\ref{sing}), we now have to calculate
$E(\Psi(\th)(1-\ga^2\cos^2(\th/2))^{-1/2})$. 
(Note that $E(C\nu(\th))=E(\nu(\th))$, where $C$ is any nonzero
constant.) 
Using (\ref{int}) and performing summation of geometric series, 
we obtain:
\be\eqalign{
\sum_{k=1}^\infty k \left(\ln
\{\Psi(\th)(1-\ga^2\cos^2{\th\over2})^{-1/2}\}\right)_k^2=
\sum_{k=1}^\infty k (\ln\Psi(\th))_k^2+\\
{1\over2\pi}\int_0^{2\pi}\ln\Psi(\th) {\ga^{-2}(1-\sin{\al\over2})^2
e^{-i\th} \over 1-\ga^{-2}(1-\sin{\al\over2})^2e^{-i\th}}d\th-
{1\over 4}\ln\left\{1-\ga^{-4}\left(1-\sin{\al\over2}\right)^4
\right\}.}\la{long}
\ee
The imaginary part of the integral here vanishes while the real part
can be written as
\[
-{1\over 2}(\ln\Psi(\th))_0+\ln\De(\Psi),
\]
where $\De(\cdot)$ is defined in (\ref{De1}).

Exponentiating (\ref{long}) and multipliyng the result with
$(2\sin(\al/2)/(1+\sin(\al/2))^{1/2}$, we finish calculation of $E$:
\be
\sqrt{2\sin(\al/2)\over 1+\sin(\al/2)}
E\left({\Psi(\th)\over\sqrt{1-\ga^2\cos^2(\th/2)}}\right)=
E(\Psi){\De(\Psi)\over\sqrt{H(\Psi)}}\sin^{1/4}{\al\over2}.
\ee

Separating the contribution of the point $\th_1=0$, 
we can now rewrite (\ref{sing}) as follows:
\be\eqalign{
D_{n-1}(\hat F)\sim
\left({\ga^{2\sum_r\al_r}H(\Psi)\over 1+\sin(\al/2)}\right)^n
n^{1/4+2\al_1+2\al_1^2+\sum_1^R\al_r^2}
E(\Psi(\th)){\De(\Psi)\over\sqrt{\Psi(0)}}
\sin^{1/4}{\al\over2}
{G(3/2+2\al_1)^2\over G(2+4\al_1)}\times
\\
\prod_1^R
\left({H(\Psi)\over
    \Psi(\th_r)}{2\sin(\vp_r/2)\over 1+\sin(\al/2)}
\right)^{\al_r}
\prod_{r=2\atop r\neq R}^{R-1}
{G(\al_r+1)^2\over G(2\al_r+1)}
\prod_{r=2\atop r\neq R}^{R-1}\left|2\sin{\th_r\over2}
\right|^{-\al_r-4\al_1\al_r}
\prod_{r\neq s\neq {1,R}}\left|2\sin{\th_r-\th_s \over 2}
\right|^{-\al_r\al_s}.}\la{sing2}
\ee

From (\ref{tf}), (\ref{fF}) we obtain
\[
\tilde f(\vp)=2^{4\al_1}\ga^{2\sum_r\al_r}\Psi(\th)\prod_{r=2}^{R-1}
\left|2\sin{\th-\th_r \over 2}\right|^{2\al_r}.
\]

The ratio (\ref{ratio}) for this $\tilde f(\vp)$ can be
written as
\be
{D_{n-1}(f)\over D_{n-1}(\hat F)}
\sim
\ga^{n^2}\frac{\ga^{\sum_r\al_r}\sqrt{\Psi(0)}
\prod_{r=2}^{R-1}(2\sin(\th_r/2))^{\al_r}}
{\De(\Psi)K\sqrt{\pi n\sin(\al/2)}}
\left(1+\sin{\al\over2}\right)^n
{2^{4\al_1}\Gamma(1+2\al_1)\over n^{2\al_1}}
\la{ratio2},
\ee
where
\be
\ln K=
{\sin(\al/2) \over 4\pi}\sum_{r=1}^R 2\al_r
\int_0^{2\pi}\frac{\ln |2\ga\sin((\th-\th_r)/2)|}
{1-\ga^2\cos^2(\th/2)}d\th.
\ee
Using the symmetry of $\th_r$'s we regroup the addends in the above
expression and write
\[
\ln\left|4\sin{\th-\th_r\over2}
\sin{\th+\th_r\over2}\right|
=
\ln|1-2b\cos\th|-\ln b,\qquad 1/b=2\cos\th_r.
\]
Now the integrals can be computed as before, and we get
\be
K=\prod_{r=1}^R\left({2\ga\sin(\vp_r/2)\over 1+\sin(\al/2)}
\right)^{\al_r}.
\ee
Putting this together with (\ref{ratio2}) and (\ref{sing2}), 
we obtain
\be
\eqalign{
D_{n-1}(f)\sim \ga^{n^2+2n\sum_1^R\al_r}
H(\Psi)^n n^{2\al_1^2+\sum_1^R \al_r^2}
\left(n\sin{\al\over 2}\right)^{-1/4}
{G(3/2+2\al_1)^2\over\sqrt{\pi}G(2+4\al_1)}E(\Psi)\times\\
2^{4\al_1}\Gamma(1+2\al_1)
\prod_1^R\left({H(\Psi)\over
    \psi(\vp_r)}\right)^{\al_r} 
\prod_2^{R-1}
{G(\al_r+1)^2\over G(2\al_r+1)}
\left(2\sin(\th_r/2)\right)^{-4\al_1\al_r}
\prod_{r\neq s\neq {1,R}}\left|2\sin{\th_r-\th_s \over 2}
\right|^{-\al_r\al_s}.}
\ee
Finally, a simple but cumbersome calculation which uses the symmetry of 
$\th_r$'s (cf. (\ref{fF})) gives
\be\eqalign{
\prod_{r\neq s\neq {1,R}}\left|2\sin{\th_r-\th_s \over 2}
\right|^{-\al_r\al_s}=\\
\ga^{\sum_2^{R-1}\al_r\al_s-\sum_2^{R-1}\al_r^2}
\prod_{r\neq s\neq 1,R}\left|2\sin{\vp_r-\vp_s \over 2}\right|^{-\al_r\al_s}
\prod_{r=2\atop r\neq R}^{[R/2]}
\left({\sin^2 (\vp_r/2) \over 1-\ga^{-2}\cos^2(\vp_r/2)}
\right)^{\al_r^2}.}
\ee 
Substituting this into the previous formula completes the derivation
of (\ref{arcsing}).

\end{document}